\documentclass[12pt,twoside]{article}
\usepackage{amssymb}
\font\teneufm=eufm10 scaled \magstep1
\font\seveneufm=eufm7 scaled \magstep1
\font\fiveeufm=eufm5  scaled \magstep1
\newfam\eufmfam
\textfont\eufmfam=\teneufm
\scriptfont\eufmfam=\seveneufm
\scriptscriptfont\eufmfam=\fiveeufm
\def\frak#1{{\fam\eufmfam\relax#1}}

\usepackage{mathrsfs}

\newfam\msbfam
\font\tenmsb=msbm10 scaled \magstep1  \textfont\msbfam=\tenmsb
\font\sevenmsb=msbm7 scaled \magstep1 \scriptfont\msbfam=\sevenmsb
\font\fivemsb=msbm5 scaled \magstep1  \scriptscriptfont\msbfam=\fivemsb
\def\Bbb{\fam\msbfam \tenmsb}

\def\RR{{\Bbb R}}
\def\CC{{\Bbb C}}
\def\BB{{\Bbb B}}

\def\NN{{\Bbb N}}

\def\PP{{\Bbb P}}
\def\HH{{\Bbb H}}

\def\ra{\rightarrow}

 \def\HollowBoxx #1#2#3{{\dimen0=#1 \advance\dimen0 by -#2
       \dimen1=#1 \advance\dimen1 by #3
        \vrule height 0pt depth #3 width #2
       \hskip -#3
       \vrule height #1 depth #3 width #3}}
 \def\LeftContraction{\mathord{\kern1.45pt \HollowBoxx{6pt}{3.5pt}{.4pt}}\,}

 \def\HollowBox #1#2#3{{\dimen0=#1 \advance\dimen0 by -#3
       \dimen1=#1 \advance\dimen1 by #3
        \vrule height #1 depth #3 width #3
        \vrule height 0pt depth #3 width #2
        \hskip -#3}}
 \def\RightContraction{\mathord{\, \HollowBox{6pt}{3.1pt}{.4pt}} \kern1.6pt}

\def\qed{{\hfill $\Box$}}
\newtheorem{theorem}{THEOREM}[section]

\newtheorem{remark}[theorem]{Remark}

\begin{document}

\begin{center}
{\Large \bf Proper Actions of High-Dimensional Groups
\medskip\\
on Complex Manifolds}
\medskip \\
\normalsize A. V. Isaev
\end{center}

\begin{quotation} \small \sl We explicitly classify all pairs $(M,G)$, where $M$ is a connected complex manifold of dimension $n\ge 2$ and $G$ is a connected Lie group acting properly and effectively on $M$ by holomorphic transformations and having dimension $d_G$ satisfying $n^2+2\le d_G<n^2+2n$. These results extend -- in the complex case -- the classical description of manifolds admitting proper actions of groups of sufficiently high dimensions. They also generalize some of the author's earlier work on Kobayashi-hyperbolic manifolds with high-dimensional holomorphic automorphism group. 
\end{quotation}

\thispagestyle{empty}

\pagestyle{myheadings}
\markboth{A. V. Isaev}{Proper Actions on Complex Manifolds}

\setcounter{section}{-1}

\section{Introduction}
\setcounter{equation}{0}

Let $M$ be a connected $C^{\infty}$-smooth manifold and $\hbox{Diff}(M)$ the group of $C^{\infty}$-smooth diffeomorphisms of $M$ endowed with the compact-open topology. A topological group $G$ is said to act continuously on $M$ by diffeomorphisms, if a continuous homomorphism $\Phi: G\ra\hbox{Diff}(M)$ is specified. The continuity of $\Phi$ is equivalent to the continuity of the action map
$$
\hat\Phi:\,G\times M\ra M,\quad (g,p)\mapsto \Phi(g)(p)=:gp.
$$ 
We only consider effective actions, that is, assume that the kernel of $\Phi$ is trivial.

The action of $G$ on $M$ is called {\it proper}, if the map
$$
\Psi:\, G\times M\ra M\times M, \quad (g,p)\mapsto (gp,p),
$$
is proper, i.e. for every compact subset $C\subset M\times M$ its inverse image $\Psi^{-1}(C)\subset G\times M$ is compact as well. For example, the action is proper if $G$ is compact. The properness of the action implies that: (i) $G$ is locally compact, hence by \cite{BM1}, \cite{BM2} (see also \cite{MZ}) it carries the structure of a Lie group and the action map $\hat\Phi$ is smooth; (ii) $\Phi$ is a topological group isomorphism between $G$ and $\Phi(G)$; (iii) $\Phi(G)$ is a closed subgroup of $\hbox{Diff}(M)$ (see \cite{Bi} for a brief survey on proper actions). Thus, one can assume that $G$ is a Lie group acting smoothly and properly on the manifold $M$, and that it is realized as a closed subgroup of $\hbox{Diff}(M)$.

Suppose now that $M$ is equipped with a Riemannian metric ${\cal G}$, and let $\hbox{Isom}(M,{\cal G})$ be the group of all isometries of $M$ with respect to ${\cal G}$. It was shown in \cite{MS} that $\hbox{Isom}(M,{\cal G})$ acts properly on $M$ (and so does its every closed subgroup). Conversely, by \cite{Pal} (see also \cite{Al}), for any Lie group acting properly on $M$ there exists a $C^{\infty}$-smooth $G$-invariant metric ${\cal G}$ on $M$. It then follows that Lie groups acting properly and effectively on the manifold $M$ by diffeomorphisms are precisely closed subgroups of $\hbox{Isom}(M,{\cal G})$ for all possible smooth Riemannian metrics ${\cal G}$ on $M$.

If $G$ acts properly on $M$, then for every $p\in M$ its isotropy subgroup
$$
G_p:=\left\{g\in G: gp=p\right\}
$$
is compact in $G$. Then by \cite{Bo} the isotropy representation
\begin{equation}
\alpha_p:\, G_p\ra GL(\RR,T_p(M)),\quad g\mapsto dg(p)\label{isotreprs}
\end{equation}
is continuous and faithful, where $T_p(M)$ denotes the tangent space to $M$ at $p$ and $dg(p)$ is the differential of $g$ at $p$. In particular, the linear isotropy subgroup
$$
LG_p:=\alpha_p(G_p)
$$
is a compact subgroup of $GL(\RR,T_p(M))$ isomorphic to $G_p$. In some coordinates in $T_p(M)$ the group $LG_p$ becomes a subgroup of the orthogonal group $O_m(\RR)$, where $m:=\hbox{dim}\,M$. Hence $\hbox{dim}\,G_p\le\hbox{dim}\, O_m(\RR)=m(m-1)/2$. Furthermore, for every $p\in M$ its orbit
$$
Gp:=\left\{gp:g\in G\right\}
$$
is a closed submanifold of $M$, and $\hbox{dim}\, Gp\le m$. Thus, setting $d_G:=\hbox{dim}\,G$, we obtain
$$
d_G=\hbox{dim}\, G_p+\hbox{dim}\,Gp\le m(m+1)/2.
$$

It is a classical result (see \cite{F}, \cite{C}, \cite{Ei}) that if $G$ acts properly on a smooth manifold $M$ of dimension $m\ge 2$ and $d_G=m(m+1)/2$, then $M$ is isometric (with respect to some $G$-invariant metric) either to one of the standard complete simply-connected spaces of constant sectional curvature $\RR^m$, $S^m$, $\HH^m$ (where $\HH^m$ is the hyperbolic space), or to $\RR\PP^m$. Next, it was shown in \cite{Wa} (see also \cite{Eg}, \cite{Y1}) that a group $G$ with $m(m-1)/2+1<d_G<m(m+1)/2$ cannot act properly on a smooth manifold $M$ of dimension $m\ne 4$. The exceptional 4-dimensional case was considered in \cite{Ish}; it turned out that a group of dimension 9 cannot act properly on a 4-dimensional manifold, and that if a 4-dimensional manifold admits a proper action of an 8-dimensional group $G$, then it has a $G$-invariant complex structure. Invariant complex structures will be discussed below in detail.

There exists also an explicit classification of pairs $(M,G)$, where $m\ge 4$, $G$ is connected, and  $d_G=m(m-1)/2+1$ (see \cite{Y1}, \cite{Ku}, \cite{O}, \cite{Ish}). Further, in \cite{KN} a reasonably explicit classification of pairs $(M,G)$, where $m\ge 6$, $G$ is connected, and  $(m-1)(m-2)/2+2\le d_G\le m(m-1)/2$, was given. We also mention a classification of $G$-homogeneous manifolds for $m=4$, $d_G=6$ (see \cite{Ish}) and a classifications of $G$-homogeneous simply-connected manifolds in the cases $m=3$, $d_G=3,4$ and $m=4$, $d_G=5$ (see \cite{C}, \cite{Pat}) obtained by E. Cartan's method of adapted frames introduced in \cite{C}. There are many other results, especially for compact subgroups, but -- to the best of our knowledge -- no complete classifications exist beyond dimension $(m-1)(m-2)/2+2$ (see \cite{Ko2}, \cite{Y2} and references therein for further details).

We study proper group actions in the complex setting with the general aim to build a theory for group dimensions lower than $(m-1)(m-2)/2+2$, thus extending -- in this setting -- the classical results mentioned above. In our setting real Lie groups act by holomorphic transformations on complex manifolds. Thus, from now on, $M$ will denote a complex manifold of complex dimension $n$ (hence $m=2n$)  and $G$ will be a subgroup of $\hbox{Aut}(M)$, the group of all holomorphic automorphisms of $M$. We will be classifying pairs $(M,G)$, but we will not be concerned with determining $G$-invariant Riemannian metrics on $M$. 

Proper actions by holomorphic transformations are found in abundance. A fundamental result due to Kaup (see \cite{Ka}) states that every closed subgroup of $\hbox{Aut}(M)$ that preserves a continuous distance on $M$ acts properly on $M$. Thus, Lie groups acting properly and effectively on $M$ by holomorphic transformations are precisely those closed subgroups of $\hbox{Aut}(M)$ that preserve continuous distances on $M$. In particular, if $M$ is a Kobayashi-hyperbolic manifold, then $\hbox{Aut}(M)$ is a Lie group acting properly on $M$ (see also \cite{Ko1}).

In the complex setting, in some coordinates in $T_p(M)$ the group $LG_p$ becomes a subgroup of the unitary group $U_n$. Hence $\hbox{dim}\,G_p\le\hbox{dim}\, U_n=n^2$, and therefore
$$
d_G\le n^2+2n.
$$
We note that $n^2+2n<(m-1)(m-2)/2+2$ for $m=2n$ and $n\ge 5$. Thus, the group dimension range that arises in the complex case, for $n\ge 5$ lies strictly below the dimension range considered in the classical real case and therefore is not covered by the existing results. Furthermore, overlaps with these results for $n=3,4$ and $n=2$, $d_G=6$ occur only in relatively easy situations and do not lead to any significant simplifications in the complex case. The only interesting overlap with the real case occurs for $n=2$, $d_G=5$ (see \cite{Pat}), but we do not discuss it in this paper. Note that in the situations when overlaps do occur, the existing classifications in the real case do not necessarily immediately lead to classifications in the complex case, since the determination of all $G$-invariant complex structures on the corresponding real manifolds may be a non-trivial task.

It was shown by Kaup in \cite{Ka} that if $d_G=n^2+2n$, then $M$ is holomorphically equivalent (in fact, holomorphically isometric with respect to some $G$-invariant metric) to one of $\BB^n:=\left\{z\in\CC^n: |z|<1\right\}$, $\CC^n$, $\CC\PP^n$, and an equivalence map $F$ can be chosen so that the group\linebreak $F\circ G\circ F^{-1}:=\left\{F\circ g\circ F^{-1}:g\in G\right\}$ is, respectively, one of the groups $\hbox{Aut}(\BB^n)$, $G(\CC^n)$, $G(\CC\PP^n)$. Here $\hbox{Aut}(\BB^n)\simeq PSU_{n,1}:=SU_{n,1}/\hbox{(center)}$ is the group of all transformations 
$$
z\mapsto\displaystyle\frac{Az+b}{cz+d},\\
$$
where
$$
\left(\begin{array}{cc}
A& b\\
c& d
\end{array}
\right)
\in SU_{n,1};
$$
$G\left(\CC^n\right)\simeq U_n\ltimes\CC^n$ is the group of all holomorphic automorphisms of $\CC^n$ of the form
\begin{equation}
z\mapsto Uz+a,\label{groupcn}
\end{equation}
where $U\in U_n$, $a\in\CC^n$ (we usually write $G\left(\CC\right)$ instead of $G\left(\CC^1\right)$); and $G\left(\CC\PP^n\right)\simeq PSU_{n+1}:=SU_{n+1}/\hbox{(center)}$ is the group of all holomorphic automorphisms of $\CC\PP^n$ of the form
$$
\zeta\mapsto U\zeta,
$$
where $\zeta$ is a point in $\CC\PP^n$ written in homogeneous coordinates, and $U\in SU_{n+1}$ (this group is a maximal compact subgroup of the complex Lie group $\hbox{Aut}(\CC\PP^n)\simeq PSL_{n+1}(\CC):=SL_{n+1}(\CC)/\hbox{(center)}$). In the above situation we say for brevity that $F$ {\it transforms}\, $G$ into one of $\hbox{Aut}(\BB^n)$, $G(\CC^n)$, $G(\CC\PP^n)$, respectively, and, in general, if $F: M_1\ra M_2$ is a biholomorphic map, $G_j\subset\hbox{Aut}(M_j)$, $j=1,2$, are subgroups and $F\circ G_1\circ F^{-1}=G_2$, we say that $F$ transforms $G_1$ into $G_2$.

We remark that the groups $\hbox{Aut}(\BB^n)$, $G(\CC^n)$, $G(\CC\PP^n)$ are the full groups of holomorphic isometries of the Bergman metric on $\BB^n$, the flat metric on $\CC^n$, and the Fubini-Study metric on $\CC\PP^n$, respectively, and that the above result due to Kaup can be obtained directly from the classification of Hermitian symmetric spaces (cf. \cite{Ak}, pp. 49--50).    

We are interested in characterizing pairs $(M,G)$ for $d_G<n^2+2n$, where $G\subset\hbox{Aut}(M)$ acts on $M$ properly. In \cite{IKra}, \cite{I1}, \cite{I2}, \cite{I3} we considered the special case where $M$ is a Kobayashi-hyperbolic manifold and $G=\hbox{Aut}(M)$, and explicitly determined all manifolds with $n^2-1\le d_{\hbox{\small Aut}(M)}< n^2+2n$, $n\ge 2$ (see \cite{I4} for a comprehensive exposition of these results). The case $d_{\hbox{\small Aut}(M)}=n^2-2$ represents the first obstruction to the existence of an explicit classification, namely, there is no good description of hyperbolic manifolds with $n=2$, $d_{\hbox{\small Aut}(M)}=2$ (see \cite{I3}, \cite{I4}); it is possible, however, that a reasonable classification exists in this case for $n\ge 3$. Our immediate goal is to generalize these results to arbitrary proper actions on not necessarily Kobayashi-hyperbolic manifolds by classifying all pairs $(M,G)$ with $n^2-1\le d_G< n^2+2n$, $n\ge 2$, where $G$ is assumed to be connected.

This classification problem splits into two cases: that of $G$-homogeneous manifolds and that of non-$G$-homogeneous ones (note that due to \cite{Ka} $G$-homogeneity always takes place for $d_G>n^2$). While the techniques that we developed for non-homogeneous Kobayashi-hyperbolic manifolds seem to work well for general non-transitive proper actions, there is a substantial difference in the homogeneous case. Indeed, due to \cite{N} every homogeneous Kobayashi-hyperbolic manifold is holomorphically equivalent to a Siegel domain of the second kind, and therefore such manifolds can be studied by using techniques available for Siegel domains (see e.g. \cite{S}). This is how homogeneous Kobayashi-hyperbolic manifolds with $n^2-1\le d_{\hbox{\small Aut}(M)}< n^2+2n$, $n\ge 2$, were determined in \cite{IKra}, \cite{I1}, \cite{I3}, \cite{I4}. This approach cannot be applied to general transitive proper actions, and one motivation for the present work is to re-obtain the classification of homogeneous Kobayashi-hyperbolic manifolds without using the non-trivial result of \cite{N}.   

For general $G$-homogeneous manifolds we have
\begin{equation}
\hbox{dim}\, G_p=d_G-2n.\label{dimestisotrop}
\end{equation}
Hence for $n^2-1\le d_G< n^2+2n$ we have $n^2-2n-1\le\hbox{dim}\, G_p<n^2$. The starting point of our method of studying $G$-homogeneous manifolds with compact isotropy subgroups within the above dimension range is describing connected subgroups of the unitary group $U_n$ of respective dimensions, thus determining the connected identity components of all possible linear isotropy subgroups. In the present paper we deal with manifolds equipped with proper actions for which $n^2+2\le d_G< n^2+2n$. Due to \cite{Ka}, all such manifolds are $G$-homogeneous, and our proofs use the description of connected closed subgroups $H\subset U_n$ with $n^2-2n+2\le\hbox{dim}\,H<n^2$ obtained in \cite{IKra} (see also \cite{I4}).   

The first step towards a general classification for proper actions with $d_G<n^2+2n$ was in fact made in \cite{IKra} where we observed that if $d_G\ge n^2+3 $, then, as in the case $d_G=n^2+2n$, the manifold must be holomorphically equivalent to one of $\BB^n$, $\CC^n$, $\CC\PP^n$. However, in \cite{IKra} we did not investigate the question what groups (if any) are possible for each of these three manifolds within the dimension range $n^2+3\le d_G<n^2+2n$. We resolve this question in Theorem \ref{trmn2+3} (see Section \ref{n2+3}). Furthermore, in Theorem \ref{trmn2+2} we give a complete classification of all pairs $(M,G)$ with $d_G=n^2+2$ (see Section \ref{n2+2}).

In the proofs of Theorems \ref{trmn2+3} and \ref{trmn2+2} we do not use the existing structure theory for actions of Lie groups on complex manifolds (see e.g. \cite{HO}, \cite{Wo}). Neither do we use the classification of Hermitian symmetric spaces due to E. Cartan (see \cite{H}), a reference to which can significantly simplify the proof of Part (ii) of Theorem \ref{trmn2+3} and that of Theorem \ref{trmn2+2} (see Remark \ref{symmspaces}). We deliberately do not refer to these general facts and give proofs based on elementary calculations involving holomorphic fundamental vector fields of the $G$-action.

Working with lower values of $d_G$ requires, in particular, further analysis of subgroups of $U_n$. For example, for the case $d_G=n^2+1$ one needs a description of closed connected $(n-1)^2$-dimensional subgroups. A description of such subgroups was given in Lemma 2.1 of \cite{IKru}, and we will attempt to deal with the case $d_G=n^2+1$ in our future work. There are a large number of examples of actions with $d_G=n^2+1$, and at this stage it is not clear whether all such actions can be classified in a reasonable way.

\section{The case $n^2+3\le d_G<n^2+2n$}\label{n2+3}
\setcounter{equation}{0}

In this section we prove the following theorem.

\begin{theorem}\label{trmn2+3}\sl Let $M$ be a connected complex manifold of dimension $n\ge 2$ and $G\subset\hbox{Aut}(M)$ a connected Lie group that acts properly on $M$ and has dimension $d_G$ satisfying $n^2+3\le d_G<n^2+2n$. Then one of the following holds:
\vspace{0.3cm}

\noindent (i) $M$ is holomorphically equivalent to $\CC^n$ by means of a map that transforms $G$ into the group $G_1(\CC^n)$ which consists of all maps of the form (\ref{groupcn}) with $U\in SU_n$ (here $d_G=n^2+2n-1$);
\vspace{0.3cm}

\noindent (ii) $n=4$ and $M$ is holomorphically equivalent to $\CC^4$ by means of a map that transforms $G$ into the group $G_2(\CC^4)$ which consists of all maps of the form (\ref{groupcn}) for $n=4$ with $U\in e^{i\RR}Sp_2$ (here $d_G=n^2+3=19$).\footnote{Here $Sp_2$ denotes the standard compact real form of $Sp_4(\CC)$.}
\end{theorem}

\noindent {\bf Proof:} Fix $p\in M$. It follows from (\ref{dimestisotrop}) that $n^2-2n+3\le \hbox{dim}\, LG_p<n^2$. Choose coordinates in $T_p(M)$ so that $LG_p\subset U_n$. Then Lemma 2.1 in \cite{IKra} (see also Lemma 1.4 in \cite{I4}) implies that the connected identity component $LG_p^0$ of $LG_p$ either is $SU_n$ or, for $n=4$, is conjugate in $U_4$ to $e^{i\RR}Sp_2$. In both cases, it follows that $LG_p$ acts transitively on directions in $T_p(M)$, that is, for any two non-zero vectors $v_1,v_2\in T_p(M)$ there exists $h\in LG_p$ such that $hv_1=\lambda v_2$ for some $\lambda\in\RR^*$ (observe that the standard action of $Sp_2$ on $\CC^4$ is transitive on the sphere $S^7=\partial\BB^4$). Now the result of \cite{GK} gives that if $M$ is non-compact, it is holomorphically equivalent to one of $\BB^n$, $\CC^n$, and an equivalence map can be chosen so that it maps $p$ into the origin and transforms $G_p$ into a subgroup of $U_n\subset G(\CC^n)$; it then follows that one can find an equivalence map that transforms $G_p^0$ either into $SU_n$ or, for $n=4$, into $e^{i\RR}Sp_2$. Furthermore, the result of \cite{BDK} gives that if $M$ is compact, it is holomorphically equivalent to $\CC\PP^n$.

Suppose first that $LG_p^0=SU_n$. In this case $d_G=n^2+2n-1$. If $M$ is holomorphically equivalent to $\BB^n$, then the equivalence map transforms $G$ into a closed subgroup of codimension 1 in $\hbox{Aut}(\BB^n)$. However, the Lie algebra of $\hbox{Aut}(\BB^n)$ is isomorphic to ${\frak {su}}_{n,1}$, and it was shown in \cite{EaI} that for $n\ge 2$ this algebra does not have codimension 1 subalgebras. Thus, $M$ cannot be equivalent to $\BB^n$. Next, if $M$ is equivalent to $\CC\PP^n$, the group $G$ is compact.  Therefore, the equivalence map transforms $G$ into a closed codimension 1 subgroup of a maximal compact subgroup in $\hbox{Aut}(\CC\PP^n)$. It then follows that one can find an equivalence map that transforms $G$ into a closed codimension 1 subgroup of $G(\CC\PP^n)$. Since $G(\CC\PP^n)$ is isomorphic to $PSU_{n+1}$, we obtain that $SU_{n+1}$ has a closed codimension 1 subgroup, which contradicts Lemma 2.1 in \cite{IKra} (see also Lemma 1.4 in \cite{I4}). Thus, $M$ cannot be equivalent to $\CC\PP^n$ either.

Assume finally that $M$ is equivalent to $\CC^n$ and let $F$ be an equivalence map that transforms $G_p^0$ into $SU_n\subset G(\CC^n)$. We will show that $F$ transforms $G$ into $G_1(\CC^n)$. We only give a proof for $n=2$ (hence $d_G=7$); the general case follows by considering copies of $SU_2$ lying in $SU_n$, and we omit details.

Denote by $(z,w)$ coordinates in $\CC^2$ and let ${\frak g}$ be the Lie algebra of holomorphic vector fields on $\CC^2$ that are fundamental vector fields of the action of $G^F:=F\circ G\circ F^{-1}$, that is, ${\frak g}$ consists of all vector fields $X$ on $\CC^2$ for which there exists an element $a$ of the Lie algebra of $G^F$ such that for all $(z,w)\in\CC^2$ we have
$$
X(z,w)=\frac{d}{dt}\Bigl[\exp(ta)(z,w)\Bigr]\Bigr|_{t=0}.
$$
Since $G^F$ acts on $\CC^2$ transitively, the algebra ${\frak g}$ is generated by ${\frak {su}}_2$ (realized as the algebra of fundamental vector fields of the standard action of $SU_2$ on $\CC^2$), and some vector fields
$$
Y_j=f_j\,{\partial}/{\partial z}+g_j\, {\partial}/{\partial w},\quad j=1,2,3,4.
$$
Here the functions $f_j$, $g_j$, $j=1,2,3,4$, are holomorphic on $\CC^2$ and satisfy the conditions
$$
\begin{array}{ll}
f_1(0)=1, & g_1(0)=0,\\
f_2(0)=i, & g_2(0)=0,\\
f_3(0)=0, & g_3(0)=1,\\ 
f_4(0)=0, & g_4(0)=i.
\end{array}
$$
To prove that $G^F=G_1(\CC^n)$, it is sufficient to show that $Y_j$ can be chosen as follows:
\begin{equation}
\begin{array}{lll}
Y_1&=& \hspace{0.2cm}{\partial}/{\partial z},\\
Y_2&=& i\, {\partial}/{\partial z},\\
Y_3&=& \hspace{0.2cm}{\partial}/{\partial w},\\
Y_4&=& i\, {\partial}/{\partial w}.
\end{array}\label{goodforms}
\end{equation}

We fix the following generators in ${\frak{su}}_2$:
$$
\begin{array}{lllll}
X_1&:=& \hspace{0.15cm}iz\, {\partial}/{\partial z}& - &iw\, {\partial}/{\partial w},\\
X_2&:=& \hspace{0.15cm}w\, {\partial}/{\partial z}& - &\hspace{0.15cm}z\, {\partial}/{\partial w},\\  
X_3&:=& iw\, {\partial}/{\partial z}& + &iz\, {\partial}/{\partial w}.
\end{array}
$$
A straightforward calculation gives
$$
\begin{array}{lll}
[Y_1,X_1](0)&=&\hspace{0.33cm}(i,0),\\

[Y_2,X_1](0)&= &-(1,0).
\end{array}
$$
It then follows that
\begin{equation}
\begin{array}{llll}
Y_1& = & -[Y_2,X_1] & \hbox{(mod ${\frak{su}}_2$)},\\
Y_2&= &\hspace{0.33cm}[Y_1,X_1] & \hbox{(mod ${\frak{su}}_2$)},
\end{array}\label{y12}
\end{equation}
which implies
$$
Y_1=-[[Y_1,X_1],X_1]\,\,\,\,\hbox{(mod ${\frak{su}}_2$)}.
$$
This identity yields
\begin{equation}
\begin{array}{l}
\displaystyle\Biggl(z\,\frac{\partial f_1}{\partial z}-3w\,\frac{\partial f_1}{\partial w}-z^2\,\frac{\partial^2 f_1}{\partial z^2}+\\
\vspace{0cm}\\
\hspace{3cm}\displaystyle 2zw\,\frac{\partial^2 f_1}{\partial z\,\partial w}-
\displaystyle w^2\,\frac{\partial^2 f_1}{\partial w^2}\Biggr)\,\partial/\partial z+
\vspace{0cm}\\
\displaystyle\Biggl(-3z\,\frac{\partial g_1}{\partial z}+w\,\frac{\partial g_1}{\partial w}-z^2\,\frac{\partial^2 g_1}{\partial z^2}+\\
\vspace{0cm}\\
\hspace{3cm}\displaystyle 2zw\,\frac{\partial^2 g_1}{\partial z\,\partial w}-w^2\,\frac{\partial^2 g_1}{\partial w^2}\Biggr)\,\partial/\partial w=0\,\,\,\hbox{(mod ${\frak{su}}_2$)}.
\end{array}\label{main1}
\end{equation}
Representing the functions $f_1$ and $g_1$ as power series around the origin, plugging these representations into (\ref{main1}) and collecting terms of fixed orders, we obtain
$$
\begin{array}{l}
\displaystyle Y_1= \left(1+\sum_{m=1}^{\infty}\left(\alpha_mz^mw^m+\alpha_m'z^{m+1}w^{m-1}\right)\right)\,\partial/\partial z+\\
\displaystyle \hspace{1cm}\left(\sum_{m=1}^{\infty}\left(\beta_mz^mw^m+\beta_m'z^{m-1}w^{m+1}\right)\right)\,\partial/\partial w\,\,\,\hbox{(mod ${\frak{su}}_2$)},
\end{array}
$$
for some $\alpha_m,\alpha_m',\beta_m,\beta_m'\in\CC$. Adding to $Y_1$ an element of ${\frak{su}}_2$ if necessary, we can assume that $Y_1$ has no linear terms, that is
\begin{equation}
\begin{array}{l}
\displaystyle Y_1= \left(1+\sum_{m=1}^{\infty}\left(\alpha_mz^mw^m+\alpha_m'z^{m+1}w^{m-1}\right)\right)\,\partial/\partial z+\\
\displaystyle \hspace{1cm}\left(\sum_{m=1}^{\infty}\left(\beta_mz^mw^m+\beta_m'z^{m-1}w^{m+1}\right)\right)\,\partial/\partial w.
\end{array}\label{y1}
\end{equation}

Further, (\ref{y12}) gives
$$
Y_2=-[[Y_2,X_1],X_1]\,\,\,\,\hbox{(mod ${\frak{su}}_2$)},
$$ 
and the application of an analogous argument to $Y_2$ yields that $Y_2$ can be chosen to have the form
\begin{equation}
\begin{array}{l}
\displaystyle Y_2= \left(i+\sum_{m=1}^{\infty}\left(\tilde\alpha_mz^mw^m+\tilde\alpha_m'z^{m+1}w^{m-1}\right)\right)\,\partial/\partial z+\\
\displaystyle \hspace{1cm}\left(\sum_{m=1}^{\infty}\left(\tilde\beta_mz^mw^m+\tilde\beta_m'z^{m-1}w^{m+1}\right)\right)\,\partial/\partial w.
\end{array}\label{y2}
\end{equation}
for some $\tilde\alpha_m,\tilde\alpha_m',\tilde\beta_m,\tilde\beta_m'\in\CC$. Next, plugging (\ref{y1}), (\ref{y2}) into either of identities (\ref{y12}) implies
$$
\begin{array}{ll}
\tilde\alpha_m=\hspace{0.35cm}i\alpha_m,& \tilde\alpha_m'=-i\alpha_m',\\
\tilde\beta_m=-i\beta_m, & \tilde\beta_m'=\hspace{0.35cm}i\beta_m',
\end{array}
$$
for all $m\in\NN$. We also observe that considering $[Y_3,X_1]$ and $[Y_4,X_1]$ yields that $Y_3$, $Y_4$ can be chosen to have the forms
$$
\begin{array}{l}
\displaystyle Y_3= \left(\sum_{m=1}^{\infty}\left(\gamma_mz^mw^m+\gamma_m'z^{m+1}w^{m-1}\right)\right)\,\partial/\partial z+\\
\displaystyle \hspace{1cm}\left(1+\sum_{m=1}^{\infty}\left(\delta_mz^mw^m+\delta_m'z^{m-1}w^{m+1}\right)\right)\,\partial/\partial w,\\
\vspace{0cm}\\
\displaystyle Y_4= \left(\sum_{m=1}^{\infty}\left(-i\gamma_mz^mw^m+i\gamma_m'z^{m+1}w^{m-1}\right)\right)\,\partial/\partial z+\\
\displaystyle \hspace{1cm}\left(i+\sum_{m=1}^{\infty}\left(i\delta_mz^mw^m-i\delta_m'z^{m-1}w^{m+1}\right)\right)\,\partial/\partial w,
\end{array}
$$
for some $\gamma_m,\gamma_m',\delta_m,\delta_m'\in\CC$.

Further, computing $[Y_j,X_2]$ for $j=1,2,3,4$ and collecting terms of orders 2 and greater, we obtain
$$
\begin{array}{l}
\alpha_m=\alpha_m'=\beta_m=\beta_m'=\gamma_m=\gamma_m'=\delta_m=\delta_m'=0,\quad m\ge 2,\\
\alpha_1'=\beta_1,\,\beta_1'=\alpha_1,\,\gamma_1'=\delta_1,\,\delta_1'=\gamma_1.
\end{array}
$$
Next, we have
$$
[Y_1,Y_2]=-2i\beta_1(2z\,\partial/\partial z+w\,\partial/\partial w).
$$
Hence $[Y_1,Y_2]=0\,\hbox{(mod ${\frak{su}}_2$)}$, which can only hold if $\beta_1=0$. Similarly, considering $[Y_3,Y_4]$ leads to $\gamma_1=0$. Finally, we compute $[Y_1,Y_3]$, $[Y_1,Y_4]$ and see that $\alpha_1=\delta_1=0$. Thus, $Y_j$ chosen as above (that is, not having linear terms) are in fact given by (\ref{goodforms}), and we have obtained (i) of the theorem.

Suppose next that $n=4$ and $LG_p^0$ is conjugate in $U_4$ to $e^{i\RR}Sp_2$. In this case $d_G=n^2+3=19$. If $M$ is equivalent to $\CC\PP^4$, the group $G$ is compact.  Therefore, one can find an equivalence map that transforms $G$ into a closed 19-dimensional subgroup of $G(\CC\PP^4)$. Since $G(\CC\PP^4)$ is isomorphic to $PSU_5$, we obtain that $SU_5$ has a closed 19-dimensional subgroup, which contradicts Lemma 2.1 in \cite{IKra} (see also Lemma 1.4 in \cite{I4}). Thus, $M$ cannot be equivalent to $\CC\PP^4$.

Assume now that $n=4$, the manifold $M$ is equivalent to one of $\BB^4$, $\CC^4$, and let $F$ be an equivalence map that transforms $G_p^0$ into $e^{i\RR}Sp_2\subset G(\CC^4)$. We will show that $F$ transforms $G$ into $G_2(\CC^4)$. Let ${\frak g}$ be the Lie algebra of fundamental vector fields of the action of $G^F:=F\circ G\circ F^{-1}$ on one of $\BB^4$, $\CC^4$, respectively. Since $G^F$ contains the one-parameter subgroup $z\mapsto e^{it}z$, $t\in\RR$, the algebra ${\frak g}$ contains the vector field
$$
Z_0:=i\sum_{k=1}^4 z_k\,\partial/\partial z_k.
$$
Hilfssatz 4.8 of \cite{Ka} then gives that every vector field in ${\frak g}$ is polynomial and has degree at most 2. Since $G^F$ acts transitively on one of $\BB^4$, $\CC^4$, the algebra ${\frak g}$ is generated by $\langle Z_0\rangle\oplus{\frak{sp}}_2$ (where $\langle Z_0\rangle$ is the one-dimensional algebra spanned by $Z_0$ and ${\frak{sp}}_2$ denotes the Lie algebra of $Sp_2$ realized as the algebra of fundamental vector fields of the standard action of $Sp_2$ on $\CC^4$), and some vector fields
\begin{equation}
\begin{array}{lll}
V_j&=&\displaystyle\sum_{k=1}^4 f_j^k\,\partial/\partial z_k,\\
\vspace{0cm}&&\\ 
W_j&=&\displaystyle\sum_{k=1}^4 g_j^k\,\partial/\partial z_k,
\end{array}\label{vwvectfields}
\end{equation}
for $j=1,2,3,4$. Here $f_j^k$, $g_j^k$ are holomorphic polynomials of degree at most 2 such that
\begin{equation}
f_j^k(0)=\delta_j^k,\quad g_j^k(0)=i\delta_j^k.\label{vwfunctions}
\end{equation}
Considering $[Z_0,[V_j,Z_0]]$, $[Z_0,[W_j,Z_0]]$ instead of $V_j$, $W_j$ if necessary, we can assume that $V_j$, $W_j$, $j=1,2,3,4$, have no linear terms (see the proof of Satz 4.9 in \cite{Ka}). Thus, we have
$$
f_j^k=\delta_j^k+\hbox{second-order terms},\quad g_j^k=i\delta_j^k+\hbox{second-order terms}.
$$
To prove that $F$ maps $M$ onto $\CC^4$ and $G^F=G_2(\CC^4)$, we need to show that all the second-order terms identically vanish, that is,
\begin{equation}
f_j^k\equiv\delta_j^k,\quad g_j^k\equiv i\delta_j^k.\label{wts}
\end{equation} 

We introduce the following vector fields from $\langle Z_0\rangle\oplus{\frak{sp}}_2$:
\begin{equation}
\begin{array}{lll}
Z_1&:=& iz_1\,\partial/\partial z_1+iz_4\,\partial/\partial z_4,\\
Z_2&:=& iz_2\,\partial/\partial z_2-iz_4\,\partial/\partial z_4,\\
Z_3&:=& iz_3\,\partial/\partial z_3+iz_4\,\partial/\partial z_4,\\
Z_4&:=&\hspace{0.15cm}z_4\,\partial/\partial z_2-z_2\,\partial/\partial z_4,\\
Z_5&:=&iz_4\,\partial/\partial z_2+iz_2\,\partial/\partial z_4,\\
Z_6&:=&z_2\,\partial/\partial z_1-z_1\,\partial/\partial z_2+z_4\,\partial/\partial z_3-z_3\,\partial/\partial z_4,\\
Z_7&:=&iz_2\,\partial/\partial z_1+iz_1\,\partial/\partial z_2-iz_4\,\partial/\partial z_3-iz_3\,\partial/\partial z_4,\\
Z_8&:=&\hspace{0.15cm} z_3\,\partial/\partial z_1-z_1\,\partial/\partial z_3,\\
Z_9&:=& iz_3\,\partial/\partial z_1+iz_1\,\partial/\partial z_3.
\end{array}\label{vectfields}
\end{equation}
It is straightforward to see that the commutators $[V_1,Z_2]$ and $[V_1,Z_3]$ vanish at the origin and have no linear terms. Hence these commutators are equal to 0, which implies that $V_1$ has the form
\begin{equation}
\begin{array}{l}
V_1=(1+\alpha z_1^2)\,\partial/\partial z_1+\beta z_1z_2\,\partial/\partial z_2+(\gamma z_1z_3+\delta z_2z_4)\,\partial/\partial z_3+\\
\hspace{9.3cm}\varepsilon z_1z_4\,\partial/\partial z_4,
\end{array}\label{v1}
\end{equation}
for some $\alpha,\beta,\gamma,\delta,\varepsilon\in\CC$. Similarly, considering $[W_1,Z_2]$ and $[W_1,Z_3]$ gives        
\begin{equation}
\begin{array}{l}
W_1=(i+\alpha' z_1^2)\,\partial/\partial z_1+\beta' z_1z_2\,\partial/\partial z_2+(\gamma' z_1z_3+\delta' z_2z_4)\,\partial/\partial z_3+\\
\hspace{9.3cm}\varepsilon' z_1z_4\,\partial/\partial z_4,
\end{array}\label{w1}
\end{equation}
for some $\alpha',\beta',\gamma',\delta',\varepsilon'\in\CC$.

Consider the commutator $[V_1,Z_1]$. It is straightforward to see from (\ref{v1}) that $[V_1,Z_1]$ does not have linear terms and that $[V_1,Z_1]-W_1$ vanishes at the origin. Hence $[V_1,Z_1]=W_1$ which implies that in (\ref{w1}) we have
$$
\alpha'=-i\alpha,\, \beta'=-i\beta,\,\gamma'=-i\gamma,\,\delta'=-i\delta,\,\varepsilon'=-i\varepsilon.
$$
Then
$$
[V_1,W_1]=-2i\left(2\alpha z_1\,\partial/\partial z_1+\beta z_2\,\partial/\partial z_2+\gamma z_3\,\partial/\partial z_3+\varepsilon z_4\,\partial/\partial z_4\right).
$$
Therefore, $[V_1,W_1]=0\,\hbox{(mod $\langle Z_0\rangle\oplus{\frak{sp}}_2$)}$, which can only hold if
\begin{equation}
\varepsilon=2\alpha-\beta+\gamma.\label{condi1}
\end{equation}
Next, we compute
$$
[V_1,Z_4]=(\varepsilon-\beta)z_1z_4\,\partial/\partial z_2+\delta(z_2^2-z_4^2)\,\partial/\partial z_3+(\varepsilon-\beta)z_1z_2\,\partial/\partial z_4.
$$
Since $[V_1,Z_4]$ does not have linear terms and vanishes at the origin, it vanishes identically, that is, we have
\begin{equation}
\varepsilon=\beta,\,\delta=0.\label{condi2}
\end{equation}

Further, computing the commutators $[V_2,Z_1]$ and $[V_2,Z_3]$, we obtain analogously to (\ref{v1})
\begin{equation}
\begin{array}{l}
V_2=\kappa z_1z_2\,\partial/\partial z_1+(1+\lambda z_2^2)\,\partial/\partial z_2+\mu z_2z_3\,\partial/\partial z_3 +\\
\hspace{7.5cm}(\nu z_1z_3+\xi z_2z_4)\,\partial/\partial z_4,
\end{array}\label{v2}
\end{equation}
for some $\kappa,\lambda,\mu,\nu,\xi\in\CC$. In addition, it is straightforward to see that $[V_1,Z_6]$ does not have linear terms and that $[V_1,Z_6]+V_2$ vanishes at the origin. Hence we have
\begin{equation}
[V_1,Z_6]=-V_2\label{v1z5}.
\end{equation}
Formulas (\ref{v2}) and (\ref{v1z5}) imply
\begin{equation}
\beta=\alpha,\,\varepsilon=\gamma-\delta.\label{condi3}
\end{equation}
Now (\ref{condi1}), (\ref{condi2}), (\ref{condi3}) yield
$$
\alpha=\beta=\gamma=\delta=\varepsilon=0,
$$
and therefore
$$
V_1=\partial/\partial z_1,\ W_1=i\,\partial/\partial z_1.
$$  
It then follows from (\ref{v1z5}) that
$$
V_2=\partial/\partial z_2.
$$

Next, we have
\begin{equation}
\begin{array}{lll}
[V_1,Z_7]&=&\hspace{0.25cm}i\,\partial/\partial z_2,\\

[V_1,Z_8]&=&-\,\partial/\partial z_3,\\

[V_1,Z_9]&=&\hspace{0.25cm}i\,\partial/\partial z_3,\\

[V_2,Z_4]&=&-\,\partial/\partial z_4,\\

[V_2,Z_5]&=&\hspace{0.25cm}i\,\partial/\partial z_4,
\end{array}\label{finalids}
\end{equation}
hence the vector fields in the right-hand side of formulas (\ref{finalids}) lie in ${\frak g}$. Since $W_2-i\,\partial/\partial z_2$, $V_3-\partial/\partial z_3$, $W_3-i\,\partial/\partial z_3$, $V_4-\partial/\partial z_4$, $W_4-i\,\partial/\partial z_4$  have no linear terms and vanish at the origin, they vanish identically and we obtain
\begin{equation}
\begin{array}{lll}
W_2&=&i\,\partial/\partial z_2,\\
V_3&=&\hspace{0.15cm}\partial/\partial z_3,\\
W_3&=&i\,\partial/\partial z_3,\\
V_4&=&\hspace{0.15cm}\partial/\partial z_4,\\
W_4&=&i\,\partial/\partial z_4.
\end{array}\label{vwformula}
\end{equation}
Thus, (\ref{wts}) holds, and we have obtained (ii) of the theorem. 

The proof is complete.\qed

\section{The case $d_G=n^2+2$}\label{n2+2}
\setcounter{equation}{0}

In this section we obtain the following result.

\begin{theorem}\label{trmn2+2}\sl Let $M$ be a connected complex manifold of dimension $n\ge 2$ and $G\subset\hbox{Aut}(M)$ a connected Lie group that acts properly on $M$ and has dimension $d_G=n^2+2$. Then one of the following holds:
\vspace{0.3cm}

\noindent (i) $M$ is holomorphically equivalent to $M'\times M''$, where $M'$ is one of $\BB^{n-1}$, $\CC^{n-1}$, $\CC\PP^{n-1}$, and $M''$ is one of $\BB^1$, $\CC$, $\CC\PP^1$; an equivalence map can be chosen so that it transforms $G$ into $G'\times G''$, where $G'$ is one of $\hbox{Aut}(\BB^{n-1})$, $G(\CC^{n-1})$, $G(\CC\PP^{n-1})$, and $G''$ is one of $\hbox{Aut}(\BB^1)$, $G(\CC)$, $G(\CC\PP^1)$, respectively;
\vspace{0.3cm}

\noindent (ii) $n=4$ and $M$ is holomorphically equivalent to $\CC^4$ by means of a map that transforms $G$ into the group $G_3(\CC^4)$ which consists of all maps of the form (\ref{groupcn}) for $n=4$ with $U\in Sp_2$.
\end{theorem}

\noindent {\bf Proof:} Fix $p\in M$. It follows from (\ref{dimestisotrop}) that $\hbox{dim}\, LG_p=n^2-2n+2$. Choose coordinates in $T_p(M)$ so that $LG_p\subset U_n$. Then Lemma 2.1 in \cite{IKra} implies that the connected identity component $LG_p^0$ of $LG_p$ either is conjugate in $U_n$ to $U_{n-1}\times U_1$ or, for $n=4$, is conjugate in $U_4$ to $Sp_2$.

Suppose first that $LG_p^0$ is conjugate to $U_{n-1}\times U_1$. By Bochner's linearization theorem (see \cite{Bo}) there exist a $G_p$-invariant neighborhood ${\cal V}$ of $p$ in $M$, an $LG_p$-invariant neighborhood ${\cal U}$ of the origin in $T_p(M)$ and a biholomorphic map $F:{\cal V}\ra {\cal U}$, with $F(p)=0$, such that for every $g\in G_p$ the following holds in ${\cal V}$:
$$
F\circ g=\alpha_p(g)\circ F,
$$
where $\alpha_p$ is the isotropy representation at $p$ (see (\ref{isotreprs})). Let ${\frak g}_M$ be the Lie algebra of fundamental vector fields of the action of $G$ on $M$, and ${\frak g}_{\cal V}$ the Lie algebra of the restrictions of the elements of ${\frak g}_M$ to ${\cal V}$. Denote by ${\frak g}$ the Lie algebra of vector fields on ${\cal U}$ obtained by pushing forward the elements of ${\frak g}_{\cal V}$ by means of $F$. Observe that ${\frak g}_M$, ${\frak g}_{\cal V}$, ${\frak g}$ are naturally isomorphic, and we denote by $\varphi:{\frak g}_M\ra{\frak g}$ the isomorphism induced by $F$.   

Choose coordinates $(z_1,\dots,z_n)$ in $T_p(M)$ so that in these coordinates $LG_p^0$ is the group of matrices of the form
\begin{equation}
\left(
\begin{array}{ll}
A & 0\\
0 & e^{i\alpha}
\end{array}
\right),\label{lpc}
\end{equation}     
where $A\in U_{n-1}$, $\alpha\in\RR$. Since $F$ transforms $G_p^0|_{\cal V}$ into $LG_p^0|_{\cal U}$ and since $G$ acts transitively on $M$, the algebra ${\frak g}$ is generated by ${\frak u}_{n-1}\oplus {\frak u}_1$ and some vector fields
$$
\begin{array}{lll}
V_j&=&\displaystyle\sum_{k=1}^n f_j^k\,\partial/\partial z_k,\\
\vspace{0cm}&&\\ 
W_j&=&\displaystyle\sum_{k=1}^n g_j^k\,\partial/\partial z_k,
\end{array}
$$
for $j=1,\dots,n$, where $f_j^k$, $g_j^k$ are holomorphic functions on ${\cal U}$ such that
$$
f_j^k(0)=\delta_j^k,\quad g_j^k(0)=i\delta_j^k.
$$
Here ${\frak u}_{n-1}\oplus {\frak u}_1$ is realized as the algebra of vector fields on ${\cal U}$ of the form
\begin{equation}
\sum_{j=1}^{n-1}\left(a_{j\,1}z_1+\dots+a_{j\,n-1}z_{n-1}\right)\,\partial/\partial z_j+iaz_n\,\partial/\partial z_n,\label{uplusu}
\end{equation}
where
$$
\left(
\begin{array}{ccc}
a_{1\,1} & \dots & a_{1\,n-1}\\
\vdots & \vdots & \vdots\\
\hspace{0.3cm}a_{n-1\,1} & \dots & \hspace{0.3cm}a_{n-1\,n-1}
\end{array}
\right)\in {\frak u}_{n-1},
$$
and $a\in\RR$. 

Observe that ${\frak g}$ contains the vector field
$$
Z_0:=i\sum_{k=1}^n z_k\,\partial/\partial z_k.
$$
Therefore, due to Hilfssatz 4.8 of \cite{Ka}, every vector field in ${\frak g}$ is polynomial and has degree at most 2. Next, considering $[Z_0,[V_j,Z_0]]$, $[Z_0,[W_j,Z_0]]$ instead of $V_j$, $W_j$ if necessary, we can assume that $V_j$, $W_j$, $j=1,\dots,n$, have no linear terms. 

Furthermore, ${\frak g}$ contains the vector fields
$$
Z_k:=iz_k\,\partial/\partial z_k,\quad k=1,\dots,n.
$$
Since for each $j=1,\dots,n$ the commutators $[V_j,Z_k]$ and $[W_j,Z_k]$, with $k\ne j$, vanish at the origin and do not contain linear terms, they vanish identically, which gives
$$
\begin{array}{lll}
V_j&=&\displaystyle\sum_{k\ne j}\alpha_j^kz_kz_j\,\partial/\partial z_k+(1+\alpha_j^jz_j^2)\,\partial/\partial z_j,\\
\vspace{0cm}\\
W_j&=&\displaystyle\sum_{k\ne j}\beta_j^kz_kz_j\,\partial/\partial z_k+(i+\beta_j^jz_j^2)\,\partial/\partial z_j,
\end{array}
$$
for some $\alpha_j^k,\beta_j^k\in\CC$ (cf. the proof of Proposition 4.9 of \cite{Ka}). Next, $[V_j,Z_j]$ has no linear terms and $[V_j,Z_j](0)=(0,\dots,0,i,0,\dots,0)$, where $i$ occurs in the $j$th position. Hence $[V_j,Z_j]=W_j$ which implies
$$
\beta_j^k=-i\alpha_j^k,
$$
for all $j,k$. Now, for $j=1,\dots,n-1$ consider the commutator $[V_j,V_n]$. Clearly, the linear part ${\cal L}_j$ of this commutator must be an element of ${\frak u}_{n-1}\oplus{\frak u}_1$. It is straightforward to see that  
$$
{\cal L}_j=\alpha_n^jz_n\,\partial/\partial z_j-\alpha_j^nz_j\,\partial/\partial z_n,
$$  
which can lie in ${\frak u}_{n-1}\oplus{\frak u}_1$ only if $\alpha_j^n=\alpha_n^j=0$ (see (\ref{uplusu})).

Thus, we have shown that for $j=1,\dots,n-1$ the vector fields $V_j$, $W_j$ do not depend on $z_n$ and the vector fields $V_n,W_n$ do not depend on $z_j$. Accordingly, we have ${\frak g}={\frak g}_1\oplus{\frak g}_2$, where ${\frak g}_1$ is the ideal generated by ${\frak u}_{n-1}$ and $V_j$, $W_j$, for $j=1,\dots,n-1$, and ${\frak g}_2$ is the ideal generated by ${\frak u}_1$ and $V_n$, $W_n$.

Let $G_j$ be the connected normal (possibly non-closed) subgroup of $G$ with Lie algebra $\tilde{\frak g}_j:=\varphi^{-1}({\frak g}_j)\subset{\frak g}_M$ for $j=1,2$. Clearly, for each $j$ the subgroup $G_j$ contains $\alpha_p^{-1}\left(L_{j\,p}\right)\subset G_p^0$, where $L_{1\,p}\simeq U_{n-1}$ and $L_{2\,p}\simeq U_1$ are the subgroups of $LG_p^0$ given by $\alpha=0$ and $A=\hbox{id}$ in formula (\ref{lpc}), respectively. Consider the orbit $G_jp$, $j=1,2$. Clearly, for each $j$ there exists a neighborhood ${\cal W}_j$ of the identity in $G_j$ such that
$$
\begin{array}{l}
{\cal W}_1p=F^{-1}\left({\cal U}'\cap\{z_n=0\}\right),\\   
{\cal W}_2p=F^{-1}\left({\cal U}'\cap\{z_1=0,\dots,z_{n-1}=0\}\right),
\end{array}
$$
for some neighborhood ${\cal U}'\subset {\cal U}$ of the origin in $T_p(M)$. Thus, each $G_jp$ is a complex (possibly non-closed) submanifold of $M$, and the ideal $\tilde{\frak g}_j$ consists exactly of those vector fields from ${\frak g}_M$ that are tangent to $G_jp$ at some point (and hence at all points). 

Furthermore, for the isotropy subgroup $G_{j\,p}$ of the point $p$ with respect to the $G_j$-action we have $G_{j\,p}^0=\alpha_p^{-1}\left(L_{j\,p}\right)$, $j=1,2$.  Since $L_{j\,p}$ acts transitively on real directions in $T_p(G_jp)$ for $j=1,2$, by \cite{GK}, \cite{BDK} we obtain that $G_1p$ is holomorphically equivalent to one of $\BB^{n-1}$, $\CC^{n-1}$, $\CC\PP^{n-1}$ and $G_2p$ is holomorphically equivalent to one of $\BB^1$, $\CC$, $\CC\PP^1$.

We will now show that each $G_j$ is closed in $G$. We assume that $j=1$; for $j=2$ the proof is identical. Let ${\frak U}$ be a neighborhood of $0$ in ${\frak g}_M$ where the exponential map into $G$ is a diffeomorphism, and let ${\frak V}:=\exp({\frak U})$. To prove that $G_1$ is closed in $G$ it is sufficient to show that for some neighborhood ${\frak W}$ of $e\in G$, ${\frak W}\subset{\frak V}$, we have $G_1\cap {\frak W}= \exp(\tilde{\frak g}_1\cap {\frak U})\cap {\frak W}$. Assuming the opposite we obtain a sequence $\{g_j\}$ of elements of $G_1$ converging to $e$ in $G$ such that for every $j$ we have $g_j=\exp(a_j)$ with $a_j\in {\frak U}\setminus\tilde{\frak g}_1$. Observe now that there exists a neighborhood ${\cal V}'$ of $p$ in $M$ foliated by complex submanifolds holomorphically equivalent to $\BB^{n-1}$ in such a way that the leaf passing through $p$ lies in $G_1p$. Specifically, we take ${\cal V}':=F^{-1}({\cal U}')$ for a suitable neighborhood ${\cal U}'\subset {\cal U}$ of the origin in $T_p(M)$, and the leaves of the foliation are then given as $F^{-1}({\cal U}'\cap\{z_n=\hbox{const}\})$. For every $s\in{\cal V}'$ we denote by $N_s$ the leaf of the foliation passing through $s$. Observe that for every $s\in{\cal V}'$ vector fields from $\tilde{\frak g}_1$ are tangent to $N_s$ at every point. Let $p_j:=g_jp$. If $j$ is sufficiently large, we have $p_j\in{\cal V}'$. We will now show that $N_{p_j}\ne N_p$ for large $j$. 

Let ${\frak U}''\subset {\frak U}'\subset {\frak U}$ be neighborhoods of 0 in ${\frak g}_M$ such that: (a) $\exp({\frak U}'')\cdot\exp({\frak U}'')\subset\exp({\frak U}')$; (b) $\exp({\frak U}'')\cdot\exp({\frak U}')\subset\exp({\frak U})$; (c) ${\frak U}'=-{\frak U}'$;  (d) $G_{1\,p}\cap\exp({\frak U}')\subset\exp(\tilde{\frak g}_1\cap {\frak U}')$. We also assume that ${\cal V}'$ is chosen so that $N_p\subset\exp(\tilde{\frak g}_1\cap {\frak U}'')p$. Suppose that $p_j\in N_p$. Then $p_j=sp$ for some $s\in\exp(\tilde{\frak g}_1\cap {\frak U}'')$ and hence $t:=g_j^{-1}s$ is an element of $G_{1\,p}$. For large $j$ we have $g_j^{-1}\in\exp({\frak U}'')$. Condition (a) now implies that $t\in\exp({\frak U}')$ and hence by (c), (d) we have $t^{-1}\in\exp(\tilde{\frak g}_1\cap {\frak U}')$. Therefore, by (b) we obtain $g_j\in\exp(\tilde{\frak g}_1\cap {\frak U})$ which contradicts our choice of $g_j$. Thus, for large $j$ the leaves $N_{p_j}$ are distinct from $N_p$. Furthermore, they accumulate to $N_p\subset G_1p$. At the same time, since vector fields from $\tilde{\frak g}_1$ are tangent to every $N_{p_j}$, we have $N_{p_j}\subset G_1p$ for all $j$, and thus the orbit $G_1p$ accumulates to itself. Below we will show that this is in fact impossible thus obtaining a contradiction.  Clearly, we only need to consider the case when $G_1p$ is non-compact, that is, equivalent to one of $\BB^{n-1}$, $\CC^{n-1}$.

Since $G_{1\,p}^0$ acts on $G_1p$ effectively, by the result of \cite{GK},  the orbit $G_1p$ is holomorphically equivalent to one of $\BB^{n-1}$, $\CC^{n-1}$ by means of a map that takes $p$ into the origin and transforms $G_{1\,p}^0$ into $U_{n-1}\subset G(\CC^{n-1})$. Consider the set $S:=G_1p\cap G_2p$. The orbit $G_1p$ accumulates to itself, and therefore $S$ contains a point other than $p$. Note that $S$ does not contain any curve. Since $G_{1\,p}^0$ preserves each of $G_1p$, $G_2p$, it preserves $S$. However, the $G_{1\,p}^0$-orbit of every point in $G_1p$ other than $p$ is a hypersurface in $G_1p$ diffeomorphic to the sphere $S^{2n-3}$. This contradiction shows that in fact $S$ consists of $p$ alone, and hence $G_1$ is closed in $G$.

Thus, $G_j$ is closed in $G$ for $j=1,2$. Hence $G_j$ acts on $M$ properly and $G_jp$ is a closed submanifold of $M$ for each $j$. Recall that $G_1p$ is equivalent to one of $\BB^{n-1}$, $\CC^{n-1}$, $\CC\PP^{n-1}$ and $G_2p$ is equivalent to one of $\BB^1$, $\CC$, $\CC\PP^1$, and denote by $F_1$, $F_2$ the respective equivalence maps. Let $K_j\subset G_j$ be the ineffectivity kernel of the $G_j$-action on $G_jp$ for $j=1,2$. Clearly, $K_j\subset G_{j\,p}$ and, since $G_{j\,p}^0$ acts on $G_jp$ effectively, $K_j$ is a discrete normal subgroup of $G_j$ for each $j$ (in particular, $K_j$ lies in the center of $G_j$ for $j=1,2$). Since $d_{G_1}=n^2-1=(n-1)^2+2(n-1)$ and $d_{G_2}=3$, the results of \cite{Ka} yield that $F_1$ can be chosen to transform $G_1/K_1$ into one of $\hbox{Aut}(\BB^{n-1})$, $G (\CC^{n-1})$, $G(\CC\PP^{n-1})$, respectively, and $F_2$ can be chosen to transform $G_2/K_2$ into one of $\hbox{Aut}(\BB^1)$, $G(\CC)$, $G(\CC\PP^1)$, respectively, where $G_j/K_j$ is viewed as a subgroup of $\hbox{Aut}(G_jp)$ for each $j$.

We will now show that the subgroup $K_j$ is in fact trivial for each $j=1,2$. We only consider the case $j=1$ since for $j=2$ the proof is identical. Clearly, $K_1\setminus\{e\}\subset G_{1\,p}\setminus G_{1\,p}^0$, and if $K_1$ is non-trivial, the compact group $G_{1\,p}$ is disconnected. Observe that any maximal compact subgroup of each of $\hbox{Aut}(\BB^{n-1})\simeq PSU_{n-1,1}$ and $G(\CC^{n-1})\simeq U_{n-1}\ltimes\CC^{n-1}$ is isomorphic to $U_{n-1}$ and therefore, if $G_1/K_1$ is isomorphic to either of these two groups, it follows that $G_{1\,p}$ is a maximal compact subgroup of $G_1$. Since $G_1$ is connected, so is $G_{1\,p}$, and therefore $K_1$ is trivial in either of these two cases. Suppose now that $G_1/K_1$ is isomorphic to $G(\CC\PP^{n-1})\simeq PSU_n$. Then the universal cover of $G_1$ is the group $SU_n$, and let $\Pi: SU_n\ra G_1$ be a covering homomorphism. Then $\Pi^{-1}(G_{1\,p}^0)^0$ is a closed $(n-1)^2$-dimensional connected subgroup of $SU_n$. It follows from Lemma 2.1 of \cite{IKru} that $\Pi^{-1}(G_{1\,p}^0)^0$ is conjugate in $SU_n$ to the subgroup of matrices of the form
\begin{equation}
\left(
\begin{array}{cc}
1/\det B & 0\\
0 & B
\end{array}
\right),\label{matrices}
\end{equation}
where $B\in U_{n-1}$. This yields that $\Pi^{-1}(G_{1\,p}^0)^0$ contains the center of $SU_n$, hence $G_{1\,p}^0$ contains the center of $G_1$. In particular, $K_1\subset G_{1\,p}^0$ which implies that $K_1$ is trivial in this case as well. Thus, $G_1$ is isomorphic to one of $\hbox{Aut}(\BB^{n-1})$, $G(\CC^{n-1})$, $G(\CC\PP^{n-1})$ and $G_2$ is isomorphic to one of $\hbox{Aut}(\BB^1)$, $G(\CC)$, $G(\CC\PP^1)$.

We remark here that since $M$ is $G$-homogeneous and $G_j$ is normal in $G$, the discussion above remains valid for any point $q\in M$ in place of $p$; in particular, all $G_j$-orbits are pairwise holomorphically equivalent, for $j=1,2$.

Next, since ${\frak g}={\frak g}_1\oplus{\frak g}_2$, the group $G$ is a locally direct product of $G_1$ and $G_2$. We claim that $H:=G_1\cap G_2$ is trivial. Indeed, $H$ is a discrete normal subgroup of each of $G_1$, $G_2$. However, every discrete normal subgroup of each of  $\hbox{Aut}(\BB^k)$, $G (\CC^k)$, $G(\CC\PP^k)$ for $k\in\NN$ is trivial, since the center of each of these groups is trivial. Hence $H$ is trivial and therefore $G=G_1\times G_2$.

We will now show that for every $q_1,q_2\in M$ the orbits $G_1q_1$ and $G_2q_2$ intersect at exactly one point. Let $g\in G$ be an element such that $gq_2=q_1$. It can be uniquely represented in the form $g=g_1g_2$ with $g_j\in G_j$ for $j=1,2$, and therefore we have $g_2q_2=g_1^{-1}q_1$. Hence the intersection $G_1q_1\cap G_2q_2$ is non-empty. We will now prove that $G_1q\cap G_2q=\{q\}$ for every $q\in M$. Suppose that for some $q\in M$ the intersection $G_1q\cap G_2q$ contain a point $q'\ne q$. Let $g_1\in G_1$ be an element such that $g_1q=q'$. Clearly, $g_1$ preserves $G_2q$. Since $g_1\in G_1$ and $G=G_1\times G_2$, the element $g_1$ commutes with every element of $G_2$. Consider the restriction $g_1':=g_1|_{G_2q}$. Let $\hat F$ be a biholomorphic map from $G_2q$ onto one of $\BB^1$, $\CC$, $\CC\PP^1$ that transforms $G_2$ into one of $\hbox{Aut}(\BB^1)$, $G(\CC)$, $G(\CC\PP^1)$, respectively. Then $\hat F$ transforms $g_1'$ into a holomorphic automorphism of one of $\BB^1$, $\CC$, $\CC\PP^1$ that lies in the centralizer of the corresponding group. In each of the three cases we immediately see that $g_1'$ is the identity, which is a contradiction. Thus, the intersection $G_1q\cap G_2q$ consists of $q$ alone for every $q\in M$.

Let, as before, $F_1$ be a biholomorphic map from $G_1p$ onto $M'$, where $M'$ is one of $\BB^{n-1}$, $\CC^{n-1}$, $\CC\PP^{n-1}$, that transforms $G_1$ into $G'$, where $G'$ is one of $\hbox{Aut}(\BB^{n-1})$, $G(\CC^{n-1})$, $G(\CC\PP^{n-1})$, respectively, and let $F_2$ be a biholomorphic map from $G_2p$ onto $M''$, where $M''$ is one of $\BB^1$, $\CC$, $\CC\PP^1$, that transforms $G_2$ into $G''$, where $G''$ is one of $\hbox{Aut}(\BB^1)$, $G(\CC)$, $G(\CC\PP^1)$, respectively. We will now construct a biholomorphic map ${\cal F}$ from $M$ onto $M'\times M''$. For $q\in M$ consider $G_2q$ and let $r$ be the unique point of intersection of $G_1p$ and $G_2q$. Let $g\in G_1$ be an element such that $r=gp$. Then we set ${\cal F}(q):=(F_1(r), F_2(g^{-1}q))$. Clearly, ${\cal F}$ is a well-defined diffeomorphism from $M$ onto $M'\times M''$. Since the foliation of $M$ by $G_j$-orbits is holomorphic for each $j$, the map ${\cal F}$ is in fact holomorphic. By construction, ${\cal F}$ transforms $G$ into $G'\times G''$. Thus, we have obtained (i) of the theorem.          
 
Suppose now that $n=4$ and $LG_p^0$ is conjugate in $U_4$ to $Sp_2$. In this case $LG_p$ acts transitively on directions in $T_p(M)$. 
Now the result of \cite{GK} gives, as before, that if $M$ is non-compact, it is holomorphically equivalent to one of $\BB^4$, $\CC^4$, and an equivalence map can be chosen so that it maps $p$ into the origin, transforms $G_p$ into a subgroup of $U_4\subset G(\CC^4)$, and transforms $G_p^0$ into $Sp_2$. Furthermore, the result of \cite{BDK} gives, as before, that if $M$ is compact, it is holomorphically equivalent to $\CC\PP^4$.

If $M$ is equivalent to $\CC\PP^4$, arguing as in the proof of Theorem \ref{trmn2+3}, we obtain that $SU_5$ has a closed 18-dimensional subgroup. This contradicts Lemma 2.1 in \cite{IKra} (see also Lemma 1.4 in \cite{I4}), and therefore $M$ cannot be equivalent to $\CC\PP^4$.

Assume now that $n=4$, the manifold $M$ is equivalent to one of $\BB^4$, $\CC^4$ and let $F$ be an equivalence map that transforms $G_p^0$ into $Sp_2\subset G(\CC^4)$. Let ${\frak g}$ be the Lie algebra of fundamental vector fields of the action of $G^F:=F\circ G\circ F^{-1}$ on one of $\BB^4$, $\CC^4$, respectively. Since $G^F$ acts transitively on one of $\BB^4$, $\CC^4$, the algebra ${\frak g}$ is generated by ${\frak{sp}}_2$ (where, as before, ${\frak{sp}}_2$ denotes the Lie algebra of $Sp_2$ realized as the algebra of fundamental vector fields of the standard action of $Sp_2$ on $\CC^4$), and some vector fields (\ref{vwvectfields}), where $f_j^k$, $g_j^k$, $j,k=1,2,3,4$ are functions holomorphic on one of $\BB^4$, $\CC^4$, respectively, and satisfying (\ref{vwfunctions}). We will show that $F$ maps $M$ onto $\CC^4$ and transforms $G$ into $G_3(\CC^4)$. To obtain this, it is sufficient to prove that one can choose 
$$
\begin{array}{lll}
V_j&=&\hspace{0.2cm}\partial/\partial z_j,\\
W_j&=& i\, \partial/\partial z_j,
\end{array}
$$
for $j=1,2,3,4$.

In our arguments we will use the following vector fields from ${\frak{sp}}_2$: $Z_4$, $Z_5$, $Z_6$, $Z_7$, $Z_8$, $Z_9$ defined in (\ref{vectfields}), as well as the vector fields
$$
\begin{array}{lll}
Z_1'&:=& iz_2\,\partial/\partial z_2-iz_4\,\partial/\partial z_4,\\
Z_2'&:=& iz_1\,\partial/\partial z_1-iz_3\,\partial/\partial z_3
\end{array}
$$
(observe that $Z_1$, $Z_2$, $Z_3$ defined in (\ref{vectfields}) do not line in ${\frak{sp}}_2$). It is straightforward to see that $[V_1,Z_1'](0)=0$, and therefore we have
\begin{equation}
[V_1,Z_1']=0\quad\hbox{(mod ${\frak{sp}}_2$)}.\label{z1prime}
\end{equation}
Representing $f_1^k$ by a power series near the origin and denoting by $\tilde f_1^k$ the non-linear part of its expansion for $k=1,2,3,4$, from (\ref{z1prime}) we obtain
$$
\begin{array}{lll}
\tilde f_1^1&=&\displaystyle\sum_{n+2l+m\ge 2}a^1_{1\,n,l,m,l}z_1^nz_2^lz_3^mz_4^l,\\
\vspace{0cm}&&\\ 
\tilde f_1^2&=&\displaystyle\sum_{n+2l+m\ge 1}a^2_{1\,n,l+1,m,l}z_1^nz_2^{l+1}z_3^mz_4^l,\\
\vspace{0cm}&&\\ 
\tilde f_1^3&=&\displaystyle\sum_{n+2l+m\ge 2}a^3_{1\,n,l,m,l}z_1^nz_2^lz_3^mz_4^l,\\
\vspace{0cm}&&\\ 
\tilde f_1^4&=&\displaystyle\sum_{n+2l+m\ge 1}a^4_{1\,n,l,m,l+1}z_1^nz_2^lz_3^mz_4^{l+1},
\end{array}
$$
where $a_{1\,n,l,m,r}^k\in\CC$. Next, we observe
$$
\begin{array}{lll}
[V_1,Z_2'](0)&=&(i,0,0,0),\\

[W_1,Z_2'](0)&=&(-1,0,0,0).
\end{array}
$$
It then follows that
\begin{equation}
\begin{array}{llll}
V_1&=&-[W_1,Z_2'] &\hbox{(mod ${\frak{sp}}_2$)},\\
W_1&=&\hspace{0.33cm}[V_1,Z_2'] &\hbox{(mod ${\frak{sp}}_2$)},
\end{array}\label{v_1w_1}
\end{equation}
which yields
\begin{equation}
V_1=-[[V_1,Z_2'],Z_2'] \quad \hbox{(mod ${\frak{sp}}_2$)}.\label{z2prime}
\end{equation}
Formula (\ref{z2prime}) implies that $\tilde f_1^k$, $k=1,2,3,4$, in fact have the forms
$$
\begin{array}{lll}
\tilde f_1^1&=&\displaystyle\sum_{n+l\ge 1}a^1_{1\,n,l,n,l}z_1^nz_2^lz_3^nz_4^l+\sum_{n,l\ge 0}a^1_{1\,n+2,l,n,l}z_1^{n+2}z_2^lz_3^nz_4^l,\\
\vspace{0cm}&&\\ 
\tilde f_1^2&=&\displaystyle\sum_{n,l\ge 0}a^2_{1\,n,l+1,n+1,l}z_1^nz_2^{l+1}z_3^{n+1}z_4^l+\sum_{n,l\ge 0}a^2_{1\,n+1,l+1,n,l}z_1^{n+1}z_2^{l+1}z_3^nz_4^l,\\
\vspace{0cm}&&\\ 
\tilde f_1^3&=&\displaystyle\sum_{n+l\ge 1}a^3_{1\,n,l,n,l}z_1^nz_2^lz_3^nz_4^l+\sum_{n,l\ge 0}a^3_{1\,n,l,n+2,l}z_1^nz_2^lz_3^{n+2}z_4^l,\\
\vspace{0cm}&&\\ 
\tilde f_1^4&=&\displaystyle\sum_{n,l\ge 0}a^4_{1\,n,l,n+1,l+1}z_1^nz_2^lz_3^{n+1}z_4^{l+1}+\sum_{n,l\ge 0}a^4_{1\,n+1,l,n,l+1}z_1^{n+1}z_2^lz_3^nz_4^{l+1}.
\end{array}
$$
In addition, (\ref{z1prime}) and (\ref{z2prime}) imply that the linear part of $V_1$ is an element of ${\frak{sp}}_2$. Subtracting this element from $V_1$, we can assume that $V_1$ has no linear part.

Next, we consider $[V_1,Z_4]$. It is easy to see that $[V_1,Z_4](0)=0$, which yields
\begin{equation}
[V_1,Z_4]=0 \quad\hbox{(mod ${\frak{sp}}_2$)}.\label{z4}
\end{equation}
It follows from (\ref{z4}) that the forms of $\tilde f_1^k$, $k=1,2,3,4$, further simplify as
$$
\begin{array}{lll}
\tilde f_1^1&=&\displaystyle\sum_{n\ge 1}a^1_{1\,n,0,n,0}z_1^nz_3^n+\sum_{n\ge 0}a^1_{1\,n+2,0,n,0}z_1^{n+2}z_3^n,\\
\vspace{0cm}&&\\ 
\tilde f_1^2&=&\displaystyle\sum_{n\ge 0}a^2_{1\,n,1,n+1,0}z_1^nz_2z_3^{n+1}+\sum_{n\ge 0}a^2_{1\,n+1,1,n,0}z_1^{n+1}z_2z_3^n,\\
\vspace{0cm}&&\\ 
\tilde f_1^3&=&\displaystyle\sum_{n\ge 1}a^3_{1\,n,0,n,0}z_1^nz_3^n+\sum_{n\ge 0}a^3_{1\,n,0,n+2,0}z_1^nz_3^{n+2},\\
\vspace{0cm}&&\\ 
\tilde f_1^4&=&\displaystyle\sum_{n\ge 0}a^4_{1\,n,0,n+1,1}z_1^nz_3^{n+1}z_4+\sum_{n\ge 0}a^4_{1\,n+1,0,n,1}z_1^{n+1}z_3^nz_4.
\end{array}
$$

Applying the above arguments to $V_3$ in place of $V_1$ we obtain that the linear part of $V_3$ at the origin is an element of ${\frak{sp}}_2$ and that the non-linear parts $\tilde f_3^k$ of the expansions around the origin of the functions $f_3^k$, $k=1,2,3,4$, have the following forms 
\begin{equation}
\begin{array}{lll}
\tilde f_3^1&=&\displaystyle\sum_{n\ge 1}a^1_{3\,n,0,n,0}z_1^nz_3^n+\sum_{n\ge 0}a^1_{3\,n+2,0,n,0}z_1^{n+2}z_3^n,\\
\vspace{0cm}&&\\ 
\tilde f_3^2&=&\displaystyle\sum_{n\ge 0}a^2_{3\,n+1,1,n,0}z_1^{n+1}z_2z_3^n+\sum_{n\ge 0}a^2_{3\,n,1,n+1,0}z_1^nz_2z_3^{n+1},\\
\vspace{0cm}&&\\
\tilde f_3^3&=&\displaystyle\sum_{n\ge 1}a^3_{3\,n,0,n,0}z_1^nz_3^n+\sum_{n\ge 0}a^3_{3\,n,0,n+2,0}z_1^nz_3^{n+2},\\
\vspace{0cm}&&\\
\tilde f_3^4&=&\displaystyle\sum_{n\ge 0}a^4_{3\,n+1,0,n,1}z_1^{n+1}z_3^nz_4+\sum_{n\ge 0}a^4_{3\,n,0,n+1,1}z_1^nz_3^{n+1}z_4,
\end{array}\label{f3k}
\end{equation}
where $a_{3\,n,l,m,r}^k\in\CC$. Next, we observe
$$
[V_1,Z_8](0)=(0,0,-1,0),
$$
which gives
\begin{equation}
V_3=-[V_1,Z_8] \quad\hbox{(mod ${\frak{sp}}_2$)}.\label{z8}
\end{equation}
Formulas (\ref{f3k}) and (\ref{z8}) imply
\begin{equation}
\begin{array}{lll}
V_1&=&(1+\alpha z_1^2+az_1z_3)\,\partial/\partial z_1+(\beta z_1z_2+bz_2z_3)\,\partial/\partial z_2+\\
&&(\alpha z_1z_3+az_3^2)\,\partial/\partial z_3+(\varepsilon z_1z_4+cz_3z_4)\,\partial/\partial z_4,
\end{array}\label{exp1}
\end{equation}
for some $a,b,c,\alpha,\beta,\varepsilon\in\CC$ (cf. (\ref{v1})). Plugging this expression into (\ref{z4}) yields
\begin{equation}
\varepsilon=\beta,\quad c=b.\label{varepsbeta1}
\end{equation}

Further, if in the above argument we replace identity (\ref{z2prime}) by the identity
$$
W_1=-[[W_1,Z_2'],Z_2'] \quad \hbox{(mod ${\frak{sp}}_2$)}
$$
(which is also a consequence of (\ref{v_1w_1})) and consider $W_3$ instead of $V_3$, we obtain that $W_1$ can be chosen to have the form
\begin{equation}
\begin{array}{lll}
W_1&=&(i+\alpha' z_1^2+a'z_1z_3)\,\partial/\partial z_1+(\beta' z_1z_2+b'z_2z_3)\,\partial/\partial z_2+\\
&&(\alpha' z_1z_3+a'z_3^2)\,\partial/\partial z_3+(\varepsilon' z_1z_4+c'z_3z_4)\,\partial/\partial z_4.
\end{array}\label{exp2}
\end{equation}
for some $a',b',c',\alpha',\beta',\varepsilon'\in\CC$ (cf. (\ref{w1})). Plugging (\ref{exp1}), (\ref{exp2}) into either of identities (\ref{v_1w_1}) we obtain
\begin{equation}
\begin{array}{l}
\alpha'=-i\alpha,\,\beta'=-i\beta,\,\varepsilon'=-i\varepsilon,\\
a'=ia,\,b'=ib,\,c'=ic.
\end{array}\label{param}
\end{equation}
Then
$$
\begin{array}{lll}
[V_1,W_1]&=&-2i\Biggl(2\alpha z_1\,\partial/\partial z_1+\Bigl(\beta z_2+(a\beta-b\alpha)z_1z_2z_3\Bigr)\,\partial/\partial z_2+\\
&&\hspace{2cm}\alpha z_3\,\partial/\partial z_3+\Bigl(\varepsilon z_4+(a\varepsilon-c\alpha)z_1z_3z_4\Bigr)\,\partial/\partial z_4\Biggr).
\end{array}
$$
Therefore, $[V_1,W_1]=0\,\hbox{(mod ${\frak{sp}}_2$)}$, which can only hold if
$$
\alpha=0,\,\varepsilon=-\beta.
$$
Together with (\ref{varepsbeta1}) this implies
$$
\beta=\varepsilon=0,
$$ 
hence we have
\begin{equation}
V_1=(1+az_1z_3)\,\partial/\partial z_1+bz_2z_3\,\partial/\partial z_2+az_3^2\,\partial/\partial z_3+bz_3z_4\,\partial/\partial z_4.\label{finalv1}
\end{equation}

If in the above arguments we interchange $Z_1'$, $Z_2'$, as well as $Z_4$, $Z_8$, and use  $V_2$ in place of $V_1$, $W_2$ in place of $W_1$, $V_4$ in place of $V_3$, and $W_4$ in place of $W_3$, we obtain that $V_2$ can be chosen to have the form
\begin{equation}
V_2=dz_1z_4\,\partial/\partial z_1+(1+ez_2z_4)\,\partial/\partial z_2+dz_3z_4\,\partial/\partial z_3+ez_4^2\,\partial/\partial z_4,\label{explv2}
\end{equation}
for some $d,e\in\CC$. We will now consider $[V_1,Z_6]$. It is straightforward to see that
$$
[V_1,Z_6](0)=(0,-1,0,0),
$$
and therefore we have
\begin{equation}
V_2=-[V_1,Z_6]\quad\hbox{(mod ${\frak{sp}}_2$)}.\label{z6}
\end{equation}
Formulas (\ref{explv2}), (\ref{z6}) imply
\begin{equation}
a=b=d=e.\label{param2}
\end{equation}
We then compute
$$
[V_1,V_2]=a\left(z_4\,\partial/\partial z_1-z_3\,\partial/\partial z_2\right).
$$
Therefore $[V_1,V_2]=0\,\hbox{(mod ${\frak{sp}}_2$)}$, which can only hold if $a=0$. Hence it follows from (\ref{exp2}), (\ref{param}), (\ref{finalv1}), (\ref{explv2}), (\ref{param2}) that
$$
\begin{array}{lll}
V_1&=&\hspace{0.2cm}\partial/\partial z_1,\\
W_1&=&i\,\partial/\partial z_1,\\
V_2&=&\hspace{0.2cm}\partial/\partial z_2.
\end{array}
$$
Therefore identities (\ref{finalids}) hold, and we obtain
$$
\begin{array}{llll}
W_2&=&i\,\partial/\partial z_2 & \hbox{(mod ${\frak{sp}}_2$)},\\
V_3&=&\hspace{0.2cm}\partial/\partial z_3 & \hbox{(mod ${\frak{sp}}_2$)},\\
W_3&=&i\,\partial/\partial z_3 & \hbox{(mod ${\frak{sp}}_2$)},\\
V_4&=&\hspace{0.2cm}\partial/\partial z_4 & \hbox{(mod ${\frak{sp}}_2$)},\\
W_4&=&i\,\partial/\partial z_4 & \hbox{(mod ${\frak{sp}}_2$)}.
\end{array}
$$
Hence $W_2$, $V_3$, $W_3$, $V_4$, $W_4$ can be chosen as in formula (\ref{vwformula}), and we have obtained (ii) of the theorem.

The proof is now complete.\qed

\begin{remark}\label{symmspaces}\rm In the situations arising in Part (ii) of Theorem \ref{trmn2+3} and in both parts of Theorem \ref{trmn2+2} the group $LG_q^0$ contains the map $-\hbox{id}$ for every $q\in M$. Therefore, $M$ equipped with a $G$-invariant Hermitian metric becomes a Hermitian symmetric space. Then, with some extra work, Part (ii) of Theorem \ref{trmn2+3} as well as all of Theorem \ref{trmn2+2} follow from E. Cartan's classification of Hermitian symmetric spaces (see \cite{H}). The same applies to Part (i) of Theorem \ref{trmn2+3} if $n$ is even. We also remark that Part (i) of Theorem \ref{trmn2+3} for all $n$ follows from the results of \cite{Wo} (see Theorem 13.1 therein). Our proofs of Theorems \ref{trmn2+3} and \ref{trmn2+2} given above are elementary and do not refer to this general theory.
\end{remark}

{\obeylines
Department of Mathematics
The Australian National University
Canberra, ACT 0200
AUSTRALIA
E-mail: alexander.isaev@maths.anu.edu.au
}

\end{document}